\documentclass{article}[10pt]
\usepackage{amsmath, amsfonts}
\usepackage{anysize}
\usepackage{nopageno}
\newtheorem{theorem}{Theorem}
\newtheorem{lemma}{Lemma}
\newtheorem{conjecture}{Conjecture}
\newtheorem{corollary}{Corollary}
\marginsize{2in}{2in}{1.5in}{2.375in}
\begin{document}
\title{On edge graceful labelings of disjoint unions of $2r$-regular edge graceful graphs}
\author{\medskip Adrian Riskin\begin{footnote}{Corresponding author.}\end{footnote}
\enspace and Georgia Weidman\\
Department of Mathematics\\
Mary Baldwin College\\
Staunton, Virginia  24401\\
ariskin@mbc.edu\\
weidmange5073@mbc.edu}
\maketitle

\begin{abstract} We prove that if $G$ is a $2r$-regular edge graceful $(p,q)$ graph with $(r,kp)=1$
then $kG$ is edge graceful for odd $k$. We also prove that for certain 
specific classes of $2r$-regular edge graceful graphs
it is possible to drop the requirement that $(r,kp)=1$\end{abstract}

\section{Introduction and definitions}

A $(p,q)$ graph, which may include multiple edges, although not loops, 
 is \textit{edge graceful} provided that it is possible to label the edges with the numbers
$1$ through $q$ in such a way that the vertex labels induced by summing incident edges mod
$p$ are distinct.  This definition is due to Lo [5], who also proved a useful necessary condition:

\begin{theorem} A $(p,q)$ graph is edge graceful only if $p | q^{2} + q - \frac{p(p-1))}{2}$. \end{theorem}

A useful survey 
of results can be found in Gallian [1].  
The $k^{th}$ power of graph $G$, denoted by $G^{k}$, is obtained from $G$ by joining all vertices
of distance $\leq k$.  We denote the cycle with $n$ vertices by $C_{n}$.  
The \textit{disjoint union} $kG$ consists of $k$ disjoint copies of graph $G$.  Lee and Seah [4] have shown that 
$C_{n}^{k}$ is edge graceful for $k < \left \lfloor \frac{n}{2} \right \rfloor$ if and only if $n$ is odd
and  for $k \geq \left \lfloor \frac{n}{2} \right \rfloor$ if and only if $n$ is a multiple of 4 or $n$ is odd.
Lee, Seah, and Lo [see 3] have studied the edge gracefulness of the disjoint union of cycles, and in particular have proved that
$kC_{n}$ is edge graceful if and only if $k$ and $n$ are odd.  Wilson and Riskin [6] proved that the cartesian product
of any number of odd cycles is edge graceful.  In this paper, we generalize these results by proving

\begin{theorem} If $G$ is a $2r$-regular edge graceful $(p,q)$ graph with $(r,kp)=1$ then $kG$ is edge graceful
for odd $k$. \end{theorem}

\noindent Since both $C_{n}^{k}$ and cartesian products of odd cycles are
 $2r$-regular, this theorem provides many new examples of edge graceful graphs.
 
 A $2r$-regular edge graceful ($p,q$) graph $G$ is \textit{stri\ae form} if there is a 
 2-factorization $S_{1}, \dots, S_{r}$ of $G$ and an edge graceful labeling of $G$ such that
 when the labels of the $p$ edges of $S_{i}$ are reduced modulo $p$ they consist of
 1, 2, \dots , $p$.  Such a labeling and such a 2-factorization are a \textit{striation} of $G$
 and the individual 2-factors are \textit{stri\ae}.  Note that the edge graceful labelings of 
 $C_{n}^{k}$ given in [4], those of certain $(n,nk)$-multigraphs given in [2],
  and those of cartesian products of odd cycles given in [6] are striations.

\section{Results}

\begin{lemma} Let $G$ be a $2r$-regular edge graceful $(p,q)$ graph.  Then $p$ is odd.\end{lemma}

\noindent \textbf{Proof} Note that $q=rp$.  Hence $q^{2}+q-\frac{p(p-1)}{2} = 
r^{2}p^{2}+rp-\frac{p(p-1)}{2}$.  Thus $p \vert \frac{p(p-1)}{2}$.  Since $(p,p-1)=1$ this means $p$ is odd.
$\hfill \diamond$\\

\noindent The proof of the following lemma is similar, and hence we omit it:

\begin{lemma} Let $G$ be a $2r$-regular edge graceful graph.  If $kG$ is edge graceful then $k$ is odd. \end{lemma}

\noindent \textbf{Proof of Theorem 2:} Note that $p$ is odd and that $q=rp$.  Let $\ell_{i}, 1 \leq i \leq q$, 
be the edge labels in an edge graceful labeling of $G$.  Let $e_{i,j}$ be the edge in $kG$ in the $j^{th}$ copy of $G$
corresponding to edge $e_{i}$ in $G$.  Let $\ell_{i,j}$, the label of edge $e_{i,j}$, be given by 
$\ell_{i,j}=\ell_{i} + (j-1)q$ for $1 \leq j \leq k$.  Clearly the $kq$ edge labels are distinct.  Let $v_{i,j}$ be
the vertex in the $j^{th}$ copy of $G$ which corresponds to vertex $v_{i}$ in $G$.  Let
$E_{1}^{i}, \dots , E_{2r}^{i}$ be the labels of the edges incident with $v_{i}$ in $G$.  Then the induced label
of vertex $v_{i,j}$ is 

$$\ell(v_{i,j})
= \sum_{s=1}^{2r} (E_{s}^{i} + (j-1)q)$$

$$= \ell(v_{i}) + 2r(j-1)q$$

$$=\ell(v_{i}) + 2r^{2}(j-1)p$$

\noindent Now suppose that $\ell(v_{i,j}) = \ell(v_{s,t})$.  Then 

$$\ell(v_{i}) + 2r^{2}(j-1)p \equiv \ell(v_{s}) + 2r^{2}(t-1)p \pmod{kp}$$

\noindent and hence

$$\ell(v_{i}) \equiv \ell(v_{s}) \pmod{p}$$

\noindent Thus $i=s$.  From this it follows that 

$$2r^{2}(j-1)p \equiv 2r^{2}(t-1)p \pmod{kp}$$

\noindent Since $(r,kp)=1$ and $kp$ is odd, $(2r^{2}, kp)=1$, and thus 

$$(j-1)p \equiv (t-1)p \pmod{kp}$$.

\noindent Since $1 \leq j, t \leq k$, we have $(j-1)p = (t-1)p$, and thus $j=t$.  Therefore the induced vertex labels 
are distinct, and so the labeling is edge graceful. $\hfill \diamond$\\

\begin{theorem} Let $G$ be a $2r$-regular stri\ae form edge graceful $(p,q)$ 
graph.  If $k$ is odd then $kG$ is edge graceful.\end{theorem}

\noindent \textbf{Proof:} Let $S_{1}, \dots , S_{r}$ be the stri\ae \space of $G$.
Let $\ell_{i}^{j}$ be the label of edge $e_{i}$ in stria $j$.  Note that we assume 
$1 \leq \ell_{i}^{j} \leq p$.  Note also that $kG$ is a $(kp, krp)$ graph.  \\

\noindent Case 1: If $r=2t+1$ we define $\ell_{i,s}^{j}$, which is to be the label of the edge in the $s^{th}$ copy of $G$ which 
corresponds to edge $e_{i}$ in $S_{j}$ of $G$, as follows:

$$\ell_{i,1}^{j} = \ell_{i}^{j} + (j-1)kp$$

\noindent for $1 \leq j \leq t+1$ and

$$\ell_{i,1}^{j} = \ell_{i}^{j} + (3t+3-j)kp-p $$

\noindent for $t+2 \leq j \leq 2t+1$.\\

\noindent And finally 

$$\ell_{i,s}^{j} = \ell_{i,s-1}^{j} + p$$

\noindent for $1 \leq j \leq t+1$ and 

$$\ell_{i,s}^{j} = \ell_{i,s-1}^{j}-p$$

\noindent for $t+2 \leq j \leq 2t+1$.  Clearly this labeling scheme produces $rkp$ distinct edge
labels.  Furthermore, in moving from copy $s-1$ of $G$ to copy $s$, we add $p$ to $2t+2$ of the 
$4t+2$ edges incident with each vertex and $-p$ to $2t$ of them.  Thus the induced label of a given 
vertex in the $s^{th}$ copy is exactly $2p$ more than the induced label of the corresponding vertex in
the $(s-1)^{th}$ copy.  Since $G$ itself is edge graceful, the induced vertex labels of the first copy of
$G$ are distinct modulo $kp$ since they are distinct modulo $p$.  Furthermore, since $2p$ has additive order
$k$ in $\mathbb{Z}_{kp}$, all of the induced vertex labels are distinct, and the labeling is edge 
graceful. \\

\noindent Case 2: If $r=2t$ we define 

$$\ell_{i,1}^{j} = \ell_{i}^{j} + (j-1)kp$$

\noindent for $1 \leq j \leq t$ and

$$\ell_{i,1}^{j} = \ell_{i}^{j} + (3t+1-j)kp-p $$

\noindent for $t+1 \leq j \leq 2t$.\\

\noindent And finally 

$$\ell_{i,s}^{1} = \ell_{i,s-1}^{1} + 2p$$

\noindent and 

$$\ell_{i,s}^{j} = \ell_{i,s-1}^{j} + p$$

\noindent for $2 \leq j \leq t$ and 

$$\ell_{i,s}^{j} = \ell_{i,s-1}^{j}-p$$

\noindent for $t+1 \leq j \leq 2t$, where all $\ell_{i,s}^{j}$ are calculated
mod $kp$ so that $1 \leq \ell_{i,s}^{j} \leq kp$.  The proof that the labeling is edge graceful is very similar to Case 1, and
so is omitted. $\hfill \diamond$\\

\begin{corollary} For odd $n, k, kC_{n}^{r}$ is edge graceful.  Also if $G$ is the cartesian product of
an odd number of odd cycles and $k$ is odd, then $kG$ is edge graceful. \end{corollary}

\section{Speculations}

\noindent We feel sure that the following is true:

\begin{conjecture} The disjoint union of an odd number of copies of a $2r$-regular 
edge graceful graph is itself edge graceful. \end{conjecture}

\noindent One way to approach this problem might be to prove:

\begin{conjecture} Every $2r$-regular edge graceful graph which has a 
2-factorization into edge graceful factors is stri\ae form. \end{conjecture}

\noindent Note that this last is related to the $(n,nk)$-multigraph conjecture of Lee [see e.g. 2], which 
speculates that every such multigraph which is decomposable into Hamiltonian cycles is edge
graceful.

\section{Bibliography}

\begin{enumerate}

\item Gallian, J. Dynamic survey of graph labeling.  Electronic Journal of Combinatorics, DS6.  
http://www.combinatorics.org/Surveys/index.html

\item Ho, Y., Lee, S.M., and Seah, E.  On the edge-graceful $(n,nk)$-multigraphs conjecture.
J. Combin. Math. Combin. Comput. 9(1991), 141-147

\item Lee, S.M. New directions in the theory of edge graceful graphs.  Proc. 6th Caribbean Conference on
Combinatorics and Computing.  1991 pp. 216-231.

\item Lee, S.M. and Seah, E.  On edge gracefulness of k-th power cycles.  Congressus Numerantium,
71(1990), 237-242.

\item Lo, S. P.  On edge graceful labelings of graphs.  Congressus Numerantium, 50(1985), 231-241.

\item Wilson, S. and Riskin, A.  Edge graceful labellings of odd cycles and their products.  Bulletin of 
the ICA, 24(1998), 57-64.

\end{enumerate}

\end{document}